\documentclass[a4paper,11pt]{amsart}
\usepackage{graphicx}
\usepackage{mathptmx}
\usepackage{amsmath}
\usepackage{amssymb}
\usepackage{enumitem}
\usepackage{xcolor}
\usepackage{xparse}
\NewDocumentCommand{\eulerian}{omm}
 {%
  \genfrac<>{0pt}{}{#2}{#3}%
  \IfValueT{#1}{_{\!#1}}%
 }

 \newmuskip\pFqmuskip

\newcommand*\pFq[6][8]{%
  \begingroup 
  \pFqmuskip=#1mu\relax
  \mathchardef\normalcomma=\mathcode`,
  \mathcode`\,=\string"8000
  \begingroup\lccode`\~=`\,
  \lowercase{\endgroup\let~}\pFqcomma
  {}_{#2}F_{#3}{\left(\genfrac..{0pt}{}{#4}{#5}\bigg|#6\right)}%
  \endgroup
}
\newcommand{\pFqcomma}{{\normalcomma}\mskip\pFqmuskip}

\newtheorem{theorem}{Theorem}

\newtheorem{remark}[theorem]{Remark}

\begin{document}

\title[A finite sum and degenerate Eulerian polynomials]{A finite sum involving generalized falling factorial polynomials and degenerate Eulerian polynomials}

\author{Taekyun Kim}
\address{Department of Mathematics, Kwangwoon University, Seoul 139-701, Republic of Korea}
\email{tkkim@kw.ac.kr}

\author{DAE SAN KIM}
\address{Department of Mathematics, Sogang University, Seoul 121-742, Republic of Korea}
\email{dskim@sogang.ac.kr}


\begin{abstract}
The aim of this paper is twofold. Firstly, we investigate a finite sum involving the generalized falling factorial polynomials, in some special cases of which we express it in terms of the degenerate Stirling numbers of the second kind, the degenerate Bernoulli polynomials and the degenerate Frobenius-Euler polynomials. Secondly, we consider the degenerate Eulerian polynomials and deduce the generating function and a recurrence relation for them.
\end{abstract}

\maketitle

\section{Introduction and preliminaries}
The Stirling number of the second $S_{2}(n,k)$ is the number of ways to partition a set of $n$ objects into $k$ nonempty subsets (see \eqref{4}). The (signed) Stirling number of the first kind $S_{1}(n,k)$ is defined such that the number of permutations of $n$ elements having exactly $k$ cycles is the nonnegative integer $(-1)^{n-k}S_{1}(n,k)=|S_{1}(n,k)|$ (see \eqref{3}). The degenerate Stirling numbers of the second kind $S_{2,\lambda}(n,k)$ (see \eqref{13}) and of the first kind $S_{1,\lambda}(n,k)$ appear most naturally when we replace the powers of $x$ by the generalized falling factorial polynomials $(x)_{k,\lambda}$ (see \eqref{2}) in the defining equations (see \eqref{3}-\eqref{5}, \eqref{13}).\par
Eulerian numbers and polynomials have long been studied due to their number-theoretic and combinatorial importance. There are many ways of defining the Eulerian polynomials and numbers, incluing the ones given in \eqref{9}-\eqref{11}. Let $[n]=\left\{1,2,\dots,n\right\}$. Then a combinatorial interpretation of the Eulerian polynomial $A_{n}(t)$ is as follows. For each permutation $\sigma$ in the symmetric group $S_{n}$, the descent and the excedance of $\sigma$ are respectively given by
\begin{align*}
&D(\sigma)= \left\{i \in [n-1] : \sigma(i) > \sigma(i + 1)\right\}, \\
&E(\sigma)= \left\{i \in [n-1] : \sigma(i) > i\right\}.
\end{align*}
Then $A_{n}(t)$ is given by 
\begin{align*}
A_{n}(t)=\sum_{\sigma \in S_{n}} t^{d(\sigma)}=\sum_{\sigma \in S_{n}} t^{e(\sigma)},
\end{align*}
where $d(\sigma)=|D(\sigma)|, e(\sigma)=|E(\sigma)|$. Further, the Eulerian numbers are defined as the coefficients of $A_{n}(t)$. \par
The aim of this paper is twofold. Firstly, we investigate a finite sum involving the generalized falling factorial polynomials, in some special cases of which we express it in terms of the degenerate Stirling numbers of the second, the degenerate Bernoulli polynomials and the degenerate Frobenius-Euler polynomials. Secondly, we consider the degenerate Eulerian polynomials and deduce the generating function and a recurrence relation for them. \par
The outline of this paper is as follows. In Section 1, we recall the facts that are needed throughout this paper.
In Section 2, we derive the main Theorem 1 which expresses the finite sum over $\alpha \le k \le m$ of the product $(k)_{\alpha,\lambda}a_{k}$,  in terms of the degenerate Stirling numbers of the second. In the case of $a_{k}=1$, we find alternative expressions for the sum, two in terms of the degenerate Stirling numbers of the second kind and one in terms of the degenerate Bernoulli polynomials. On the other hand, in the case of $a_{k}=x^{k}$, the sum is represented in terms of the degenerate Frobenius-Euler polynomials. In Section 3, we consider a degenerate version of the Eulerian polynomials, namely the degenerate Eulerian polynomials introduced by Carlitz, find the exponential generating function and the recurrence relation of them and express the sum $\sum_{k=1}^{m}(k)_{n,\lambda}x^{k}$ in terms of the degenerate Eulerian polynomials. \par
For any $\lambda\in\mathbb{R}$, the degenerate exponentials are defined by 
\begin{equation}
e_{\lambda}^{x}(t)=\sum_{n=0}^{\infty}(x)_{n,\lambda}\frac{t^{n}}{n!},\quad (\mathrm{see}\ [11]),\label{1}
\end{equation}
where the generalized falling factorial polynomials $(x)_{n,\lambda}$ are given by
\begin{equation}
(x)_{0,\lambda}=1,\quad (x)_{n,\lambda}=x(x-\lambda)(x-2\lambda)\cdots\big(x-(n-1)\lambda\big),\quad (n\ge 1).\label{2}	
\end{equation}
Note that $\displaystyle\lim_{\lambda\rightarrow 0}e_{\lambda}^{x}(t)=e^{xt},\ \ \lim_{\lambda\rightarrow 0}(x)_{n,\lambda}=x^{n}\displaystyle$. \par 
In particular, for $\displaystyle x=1,\ e_{\lambda}(t)=e_{\lambda}^{1}(t)=\sum_{n=0}^{\infty}\frac{(1)_{n,\lambda}}{n!}t^{n}\displaystyle$. \par 
It is well known that the Stirling numbers of the first kind are defined by 
\begin{equation}
(x)_{n}=\sum_{l=0}^{n}S_{1}(n,l)x^{l},\quad (n\ge 0),\quad (\mathrm{see}\ [10]),\label{3}	
\end{equation}
where $(x)_{0}=1,\quad (x)_{n}=x(x-1)\cdots\big(x-(n-1)\big),\ (n\ge 1)$. \par 
As the inversion formula of \eqref{3}, the Stirling numbers of the second kind are defined by 
\begin{equation}
x^{n}=\sum_{l=0}^{n}S_{2}(n,l)(x)_{l},\quad (n\ge 0),\quad (\mathrm{see}\ [12,13,16-20]).\label{4}
\end{equation}
Recently, Kim-Kim introduced the degenerate Stirling numbers of the first kind given by 
\begin{equation}
(x)_{n}=\sum_{l=0}^{n}S_{1,\lambda}(n,l)(x)_{l,\lambda},\quad (n\ge 0),\quad (\mathrm{see}\ [10,15]).\label{5}
\end{equation}
In addition, the unsigned degenerate Stirling numbers of first kind are defined by 
\begin{equation}
\langle x\rangle_{n}=\sum_{l=0}^{n}{n \brack l}_{\lambda}\langle x\rangle_{l,\lambda},\quad (n\ge 0),\quad (\mathrm{see}\ [13]),\label{6}	
\end{equation}
where 
\begin{align*}
&\langle x\rangle_{0}=1,\quad \langle x\rangle_{n}=x(x+1)(x+n-1),\quad (n\ge 1),\\
&\langle x\rangle_{0,\lambda}=1,\quad \langle x\rangle_{n,\lambda}=x(x+\lambda)\cdots(x+(n-1)\lambda),\quad (n\ge 1).
\end{align*}
From \eqref{5} and \eqref{6}, we note that 
\begin{equation}
{n \brack k}_{\lambda}=(-1)^{n-k}S_{1,\lambda}(n,k),\quad (n,k\ge 0),\quad (\mathrm{see}\ [13]).\label{7}
\end{equation}
The Eulerian numbers are given explicitly by the sum 
\begin{equation}
\eulerian{n}{m}=\sum_{j=0}^{m+1}(-1)^{j}\binom{n+1}{j}(m+1-j)^{n},\,\, (n \ge 1), \,\,\eulerian{0}{0}=1,\quad (\mathrm{see}\ [4-8]),\label{8}
\end{equation}
and the Eulerian polynomials are defined by
\begin{equation}
A_{n}(t)=\sum_{m=0}^{n}\eulerian{n}{m}t^{m},\quad (n\ge 0),\quad (\mathrm{see}\ [5,6,8,9]).\label{9}
\end{equation}
Actually, $\eulerian{n}{n}=0$, for $n \ge 1$, and hence $A_{n}(t)$ has degree $n-1$, for $n \ge 1$. \par
From \eqref{9}, we note that 
\begin{equation}
\frac{A_{n}(t)}{(1-t)^{n+1}}=\sum_{j=0}^{\infty}t^{j}(j+1)^{n},\quad (n\ge 0),\quad (\mathrm{see}\ [5-9],\label{10}	
\end{equation}
and 
\begin{equation}
\frac{t-1}{t-e^{(t-1)x}}=\sum_{n=0}^{\infty}A_{n}(t)\frac{x^{n}}{n!},\quad (\mathrm{see}\ [8,9]).\label{11}
\end{equation}
In [5], Carlitz considered the degenerate Bernoulli polynomials given by 
\begin{equation}
\frac{t}{e_{\lambda}(t)-1}e_{\lambda}^{x}(t)=\sum_{n=0}^{\infty}\beta_{n,\lambda}(x)\frac{t^{n}}{n!}.\label{12}
\end{equation}
When $x=0$, $\beta_{n,\lambda}=\beta_{n,\lambda}(0)$ are called the degenerate Bernoulli numbers. \par 
Note that $\displaystyle\lim_{\lambda\rightarrow 0}\beta_{n,\lambda}(x)=B_{n}(x),\ (n\ge 0)\displaystyle$, where $B_{n}(x)$are the ordinary Bernoulli polynomials given by 
\begin{displaymath}
	\frac{t}{e^{t}-1}e^{xt}=\sum_{n=0}^{\infty}B_{n}(x)\frac{t^{n}}{n!},\quad (\mathrm{see}\ [1-21]).
\end{displaymath}
Note that 
\begin{displaymath}
	\beta_{n,\lambda}(x)=\sum_{l=0}^{n}\binom{n}{l}\beta_{l,\lambda}(x)_{n-l,\lambda},\quad (n\ge 0).
\end{displaymath}
As the inversion formula of \eqref{5}, the degenerate Stirling numbers of the second kind are defined by 
\begin{equation}
(x)_{n,\lambda}=\sum_{k=0}^{n}S_{2,\lambda}(n,k)(x)_{k},\quad (n\ge 0),\quad (\mathrm{see}\ [10]).\label{13}
\end{equation}
From \eqref{13}, we note that 
\begin{equation}
\frac{1}{k!}\Big(e_{\lambda}(t)-1\Big)^{k}=\sum_{n=k}^{\infty}S_{2,\lambda}(n,k)\frac{t^{n}}{n!},\quad (\mathrm{see}\ [10]).\label{14}	
\end{equation}
By \eqref{14}, we easily get 
\begin{equation}
S_{2,\lambda}(n,k)=\frac{1}{k!}\sum_{l=0}^{k}\binom{k}{l}(-1)^{k-l}(l)_{n,\lambda},\quad (n\ge k\ge 0),\quad (\mathrm{see}\ [10,15]).\label{15}
\end{equation}
In [14], the degenerate Frobenius-Euler polynomials are defined by 
\begin{equation}
\frac{1-u}{e_{\lambda}(t)-u}e_{\lambda}^{x}(t)=\sum_{n=0}^{\infty}H_{n,\lambda}(x|u)\frac{t^{n}}{n!},\label{16}
\end{equation}
where $u(\ne 1)\in\mathbb{C}$. \par 
Note that $\lim_{\lambda\rightarrow 0}H_{n,\lambda}(x|u)=H_{n}(x|u)$ are the ordinary Frobenius-Euler polynomials given by 
\begin{displaymath}
\frac{1-u}{e^{t}-u}e^{xt}=\sum_{n=0}^{\infty}H_{n}(x|u)\frac{t^{n}}{n!},\quad (\mathrm{see}\ [1,2,14,21]).
\end{displaymath}
When $x=0$, $H_{n,\lambda}(0|u)=H_{n,\lambda}(u)$ are called the degenerate Frobenius-Euler numbers. \par 

\section{A finite sum involving generalized falling factorial polynomials}
Here we derive the main Theorem 1 which represents the finite sum over $\alpha \le k \le m$ of the product $(k)_{\alpha,\lambda}a_{k}$,  in terms of the degenerate Stirling numbers of the second. In the case of $a_{k}=1$, we express the sum in terms of the degenerate Stirling numbers of the second kind and the degenerate Bernoulli polynomials. Also, in case of $a_{k}=x^{k}$, we represent it in terms of the degenerate Frobenius-Euler polynomials. \par
Then we have 
\begin{align}
(n)_{\alpha,\lambda}&=\sum_{k=0}^{\alpha}S_{2,\lambda}(\alpha,k)(n)_{k}=\sum_{k=0}^{\alpha}S_{2,\lambda}(\alpha,k)k!\binom{n}{k}\label{17}\\
&=\sum_{k=0}^{n}S_{2,\lambda}(\alpha,k)k!\binom{n}{k},\nonumber
\end{align}
for any $n\in\mathbb{N}$ with $n\ge \alpha$. \par
By \eqref{17}, for $m \ge \alpha$, we get 
\begin{align}
\sum_{k=\alpha}^{m}(k)_{\alpha,\lambda}a_{k}&=\sum_{k=\alpha}^{m}a_{k}\sum_{j=1}^{k}S_{2,\lambda}(\alpha,j)j!\binom{k}{j}\label{19}\\
&=\sum_{j=1}^{m}S_{2,\lambda}(\alpha,j)j!\sum_{k=j}^{m}\binom{k}{j}a_{k}
-\sum_{j=1}^{\alpha-1}S_{2,\lambda}(\alpha,j)j!\sum_{k=j}^{\alpha-1} \binom{k}{j}a_{k}.\nonumber	
\end{align}
Here one should note that $S_{2,\lambda}(\alpha,0)=0$, since $\alpha$ is a positive integer.
Therefore, we obtain the following theorem. 
\begin{theorem}
Let $\{a_{n}\}_{n=1}^{m}$ be any sequence. Then we have 
\begin{displaymath}
\sum_{k=\alpha}^{m}(k)_{\alpha,\lambda}a_{k}=\sum_{k=1}^{m}S_{2,\lambda}(\alpha,k)k!\sum_{l=k}^{m}\binom{l}{k}a_{l}-\sum_{k=1}^{\alpha-1}S_{2,\lambda}(\alpha,k)k!\sum_{l=k}^{\alpha-1}\binom{l}{k}a_{l},
\end{displaymath}	
where $m \in\mathbb{N}$ with $m \ge \alpha$. 
\end{theorem}
Let us take $a_{k}=1$ in Theorem 1. Then we have 
\begin{equation}
\sum_{k=\alpha}^{m}(k)_{\alpha,\lambda}=\sum_{k=1}^{m}S_{2,\lambda}(\alpha,k)k!\sum_{l=k}^{m}\binom{l}{k}-\sum_{k=1}^{\alpha-1}S_{2,\lambda}(\alpha,k)k!\sum_{l=k}^{\alpha-1}\binom{l}{k}.\label{20}
\end{equation}
Note that 
\begin{equation}
\sum_{l=k}^{m}\binom{l}{k}=\sum_{l=k}^{m}\bigg(\binom{l+1}{k+1}-\binom{l}{k+1}\bigg)=\binom{m+1}{k+1}.\label{21}
\end{equation}
Thus, by \eqref{20} and \eqref{21}, we obtain the following theorem. 
\begin{theorem}
For $m \in\mathbb{N}$ with $m \ge \alpha$, we have 	
\begin{displaymath}
\sum_{k=\alpha}^{m}(k)_{\alpha,\lambda}=\sum_{k=1}^{m}S_{2,\lambda}(\alpha,k)k!\binom{m+1}{k+1}-\sum_{k=1}^{\alpha-1}S_{2,\lambda}(\alpha,k)k!\binom{\alpha}{k+1}.
\end{displaymath}
\end{theorem}
In [14], the degenerate $\lambda$-binomial coefficients are defined by 
\begin{equation}
\binom{k}{n}_{\lambda}=\frac{(k)_{n,\lambda}}{n!},\quad (n,k\ge 0).\label{22}	
\end{equation}
From Theorem 2 and \eqref{22}, we note that 
\begin{equation}
\sum_{k=n}^{m}\binom{k}{n}_{\lambda}=\frac{1}{n!}\sum_{k=1}^{m}S_{2,\lambda}(n,k)k!\binom{m+1}{k+1}-\frac{1}{n!}\sum_{k=1}^{n-1}S_{2,\lambda}(n,k)k!\binom{n}{k+1}.\label{23}
\end{equation}
From \eqref{5}, we note that 
\begin{equation}
\binom{n}{m}=\frac{1}{m!}\sum_{k=0}^{m}S_{1,\lambda}(m,k)(n)_{k,\lambda}.\label{24}	
\end{equation}
By \eqref{13}, we easily get 
\begin{align}
&\sum_{k=0}^{\alpha+1}S_{2,\lambda}(\alpha+1,k)(x)_{k}=(x)_{\alpha+1,\lambda}=(x)_{\alpha,\lambda}(x-\alpha\lambda)	\label{25}\\
&=\sum_{k=0}^{\alpha}S_{2,\lambda}(\alpha,k)(x)_{k}(x-\alpha\lambda)=\sum_{k=0}^{\alpha}S_{2,\lambda}(\alpha,k)(x)_{k}(x-k+k-\alpha\lambda)\nonumber \\
&=\sum_{k=0}^{\alpha+1}\Big\{S_{2,\lambda}(\alpha,k-1)+(k-\alpha\lambda)S_{2,\lambda}(\alpha,k)\Big\}(x)_{k},\nonumber
\end{align}
because $S_{2,\lambda}(\alpha,k)=0$ if $k <0, and \ S_{2,\lambda}(\alpha,k)=0$ if $k>\alpha$. \par 
By comparing the coefficients on both sides of \eqref{25}, we get 
\begin{equation}
S_{2,\lambda}(\alpha+1,k)=S_{2,\lambda}(\alpha,k-1)+(k-\alpha\lambda)S_{2,\lambda}(\alpha,k), \quad (0 \le k \le \alpha+1).\label{26}
\end{equation}
From \eqref{26}, we note that 
\begin{align}
&\sum_{j=1}^{m}\binom{m}{j}(j-1)!\Big\{S_{2,\lambda}(\alpha+1,j)+\alpha\lambda S_{2,\lambda}(\alpha,j)\Big\}\nonumber\\
&=\sum_{j=1}^{m}\binom{m}{j}(j-1)!\Big\{jS_{2,\lambda}(\alpha,j)+S_{2,\lambda}(\alpha,j-1)\Big\}\label{26-1}\\
&=\sum_{j=1}^{m}\binom{m}{j}j!S_{2,\lambda}(\alpha,j)+\sum_{j=0}^{m-1}\binom{m}{j+1}j!S_{2,\lambda}(\alpha,j)\nonumber\\
&=\sum_{j=1}^{m}\bigg(\binom{m}{j}+\binom{m}{j+1}\bigg)j!S_{2,\lambda}(\alpha,j)=\sum_{j=1}^{m}\binom{m+1}{j+1}j!S_{2,\lambda}(\alpha,j).\nonumber
\end{align}
Therefore, by Theorem 2 and \eqref{26-1}, we obtain the following theorem. 
\begin{theorem}
For $m\in\mathbb{N}$ with $m \ge \alpha$, we have 
\begin{align*}
\sum_{k=\alpha}^{m}(k)_{\alpha,\lambda}&=\sum_{k=1}^{m}\binom{m}{k}(k-1)!\Big\{S_{2,\lambda}(\alpha+1,k)+\alpha\lambda S_{2,\lambda}(\alpha,k)\Big\}\\
&-\sum_{k=1}^{\alpha-1}\binom{\alpha-1}{k}(k-1)!\Big\{S_{2,\lambda}(\alpha+1,k)+\alpha\lambda S_{2,\lambda}(\alpha,k)\Big\}.
\end{align*}	
\end{theorem}
From \eqref{12}, we note that 
\begin{align}
\frac{t}{e_{\lambda}(t)-1}e_{\lambda}^{1-x}(t)&=\frac{-t}{e_{\lambda}^{-1}(t)-1}e_{\lambda}^{-x}(t)=\frac{-t}{e_{-\lambda}(-t)-1}e_{-\lambda}^{x}(-t)\label{27} \\
&=\sum_{n=0}^{\infty}(-1)^{n}\beta_{n,-\lambda}(x)\frac{t^{n}}{n!}. \nonumber
\end{align}
Thus, by \eqref{12} and \eqref{27} we get 
\begin{displaymath}
	\beta_{n,\lambda}(1-x)=(-1)^{n}\beta_{n,-\lambda}(x),\quad (n\ge 0).
\end{displaymath}
From \eqref{12}, we have 
\begin{align}
\sum_{k=\alpha}^{m}e_{\lambda}^{k}(t)&=\frac{1}{e_{\lambda}(t)-1}\big(e_{\lambda}^{m+1}(t)-e_{\lambda}^{\alpha}(t)\big)=\frac{1}{t}\frac{t}{e_{\lambda}(t)-1}\big(e_{\lambda}^{m+1}(t)-e_{\lambda}^{\alpha}(t)\big)\label{28}\\
&=\sum_{j=0}^{\infty}\frac{1}{j+1}\big(\beta_{j+1,\lambda}(m+1)-\beta_{j+1,\lambda}(\alpha)\big)\frac{t^{j}}{j!}.\nonumber	
\end{align}
On the other hand, 
\begin{equation}
\sum_{k=\alpha}^{m}e_{\lambda}^{k}(t)=\sum_{j=0}^{\infty}\bigg(\sum_{k=\alpha}^{m}(k)_{j,\lambda}\bigg)\frac{t^{j}}{j!}.\label{29}
\end{equation}
Thus, by \eqref{28} and \eqref{29}, we get 
\begin{equation}
\sum_{k=\alpha}^{m}(k)_{\alpha,\lambda}=\frac{1}{\alpha+1}\big(\beta_{\alpha+1,\lambda}(m+1)-\beta_{\alpha+1,\lambda}(\alpha)\big).\label{30}
\end{equation}
Therefore, by Theorem 3 and \eqref{30}, we obtain the following theorem. 
\begin{theorem}
For $m\in\mathbb{N}$ with $m \ge \alpha$, we have 
\begin{align*}
\frac{1}{\alpha+1}\big(\beta_{\alpha+1,\lambda}(m+1)-\beta_{\alpha+1,\lambda}(\alpha)\big)&=\sum_{k=1}^{m}\binom{m}{k}(k-1)!\Big\{S_{2,\lambda}(\alpha+1,k)+\alpha\lambda S_{2,\lambda}(\alpha,k)\Big\}\\
&-\sum_{k=1}^{\alpha-1}\binom{\alpha-1}{k}(k-1)!\Big\{S_{2,\lambda}(\alpha+1,k)+\alpha\lambda S_{2,\lambda}(\alpha,k)\Big\}.
\end{align*}	
\end{theorem}
Let us take $a_{k}=x^{k}$ in Theorem 2. Then we have 
\begin{align}
\sum_{k=\alpha}^{m}(k)_{\alpha,\lambda}x^{k}&=\sum_{k=1}^{m}S_{2,\lambda}(\alpha,k)k!\sum_{l=k}^{m}\binom{l}{k}x^{l}-\sum_{k=1}^{\alpha-1}S_{2,\lambda}(\alpha,k)k!\sum_{l=k}^{\alpha-1}\binom{l}{k}x^{l}\label{31}\\
&=\sum_{k=1}^{m}S_{2,\lambda}(\alpha,k)k!x^{k}\sum_{l=0}^{m-k}\binom{l+k}{k}x^{l}-\sum_{k=1}^{\alpha-1}S_{2,\lambda}(\alpha,k)k!x^{k}\sum_{l=0}^{\alpha-k-1}\binom{l+k}{k}x^{l},\nonumber
\end{align}
where $m\in\mathbb{N}$ with $m \ge \alpha$. \par 
From \eqref{16}, we note that 
\begin{equation}
H_{n,\lambda}(x|u)=\sum_{k=0}^{n}\binom{n}{k}H_{k,\lambda}(u)(x)_{n-k,\lambda},\quad (n\ge 0).\label{33}
\end{equation}
Now, we observe that 
\begin{align}
\sum_{k=\alpha}^{m}e_{\lambda}^{k}(t)x^{k}&=\frac{1}{xe_{\lambda}(t)-1}\big(e_{\lambda}^{m+1}(t)x^{m+1}-e_{\lambda}^{\alpha}(t)x^{\alpha}\big)\nonumber\\
&=\frac{1}{x-1}\bigg(x^{m+1}\frac{1-x^{-1}}{e_{\lambda}(t)-x^{-1}}e_{\lambda}^{m+1}(t)-x^{\alpha}\frac{1-x^{-1}}{e_{\lambda}(t)-x^{-1}}e_{\lambda}^{\alpha}(t)\bigg)\label{34} \\
&=\frac{1}{x-1}\bigg\{\sum_{n=0}^{\infty}\Big(x^{m+1}H_{n,\lambda}(m+1|x^{-1})-x^{\alpha}H_{n,\lambda}(\alpha|x^{-1})\Big)\bigg\}\frac{t^{n}}{n!}.\nonumber
\end{align}
On the other hand, 
\begin{equation}
\sum_{k=\alpha}^{m}e_{\lambda}^{k}(t)x^{k}=\sum_{n=0}^{\infty}\bigg(\sum_{k=\alpha}^{m}(k)_{n,\lambda}x^{k}\bigg)\frac{t^{n}}{n!}. \label{35}	
\end{equation}
Therefore, by \eqref{31}, \eqref{34} and \eqref{35}, we obtain the following theorem. 
\begin{theorem}
Let $m\in\mathbb{N}$ with $m \ge \alpha$. Then we have 
\begin{displaymath}
\sum_{k=\alpha}^{m}(k)_{\alpha,\lambda}x^{k}=\frac{1}{x-1}\Big\{x^{m+1}H_{\alpha,\lambda}(m+1|x^{-1})-x^{\alpha}H_{\alpha,\lambda}(\alpha|x^{-1})\Big\}.
\end{displaymath}
In addition, we have
\begin{align*}
\frac{1}{x-1}\Big\{x^{m+1}H_{\alpha,\lambda}(m+1|x^{-1})-x^{\alpha}H_{\alpha,\lambda}(\alpha|x^{-1})\Big\}&=\sum_{k=1}^{m}S_{2,\lambda}(\alpha,k)k!x^{k}\sum_{l=0}^{m-k}\binom{l+k}{k}x^{l}\\
&-\sum_{k=1}^{\alpha-1}S_{2,\lambda}(\alpha,k)k!x^{k}\sum_{l=0}^{\alpha-k-1}\binom{l+k}{k}x^{l}.
\end{align*}
\end{theorem}

\section{Degenerate Eulerian polynomials and numbers}
In this section, we consider the degenerate Eulerian polynomials, find the exponential generating function and the recurrence relation of them and express the sum $\sum_{k=1}^{m}(k)_{n,\lambda}x^{k}$ in terms of the degenerate Eulerian polynomials. \par
We recall that the degenerate Eulerian polynomials are given by 
\begin{equation}
\frac{A_{n,\lambda}(x)}{(1-x)^{n+1}}=\sum_{j=0}^{\infty}(j+1)_{n,\lambda}x^{j},\quad(n \ge 0).\label{36} 
\end{equation}
Then $A_{0,\lambda}(x)=1$. Let $n \ge 1$. From \eqref{36}, we observe that 
\begin{align}
A_{n,\lambda}(x)&=\sum_{j=0}^{\infty}(j+1)_{n,\lambda}x^{j} \sum_{k=0}^{n+1}(-1)^{k}\binom{n+1}{k}x^{k}\label{37}\\
&=\sum_{m=0}^{\infty}\bigg(\sum_{k=0}^{m+1}(-1)^{k}\binom{n+1}{k}(m-k+1)_{n,\lambda}\bigg)x^{m}\nonumber \\ &=\sum_{m=0}^{n}\bigg(\sum_{k=0}^{m+1}(-1)^{k}\binom{n+1}{k}(m-k+1)_{n,\lambda}\bigg)x^{m}\nonumber 
\end{align}
Now, we define the degenerate Eulerian numbers as
\begin{align}
\eulerian{n}{m}_{\lambda}=\sum_{k=0}^{m+1}(-1)^{k}\binom{n+1}{k}(m-k+1)_{n,\lambda},\quad (n \ge 1), \quad \eulerian{0}{0}_{\lambda}=1.\label{38}
\end{align}
Therefore, by \eqref{37} and \eqref{38}, we obtain the following theorem. 
\begin{theorem}
For any nonnegative integer $n$, we have 
\begin{displaymath}
A_{n,\lambda}(x)=\sum_{m=0}^{n}\eulerian{n}{m}_{\lambda}x^{m}.
\end{displaymath}	
\end{theorem}

\begin{remark}
(a) Carlitz introduced the degenerate Eulerian polynomials $E_{n,\lambda}(x)$ which are slightly different from $A_{n,\lambda}(x)$ in \eqref{36} (see [5]). Indeed, they are defined as 
\begin{align*}
\frac{E_{n,\lambda}(x)}{(1-x)^{n+1}}=\sum_{j=0}^{\infty}(j)_{n,\lambda}x^{j},\quad(n \ge 0).
\end{align*}
(b) As we see, for example, from Theorem 10, $A_{n,\lambda}(t)$ has degree $n-1$, for $n \ge 1$. Thus we have
\begin{align*}
\eulerian{n}{n}_{\lambda}=\sum_{k=0}^{n+1}(-1)^{k}\binom{n+1}{k}(n-k+1)_{n,\lambda}=0,\quad \mathrm{for\,\, any}\,\, \lambda \,\,\mathrm{and}\,\, n \ge 1.
\end{align*}
\end{remark}

Let $n$ be a positive integer. Then we observe that 
\begin{equation}
	x^{n\lambda-1}\bigg(x^{1-\lambda}\frac{d}{dx}\bigg)^{n}\frac{1}{1-x}=\sum_{j=0}^{\infty}(j+1)_{n,\lambda}x^{j}\Longleftrightarrow x^{n\lambda}\bigg(x^{1-\lambda}\frac{d}{dx}\bigg)^{n}\frac{1}{1-x}=x\sum_{j=0}^{\infty}(j+1)_{n,\lambda}x^{j},\label{39}
\end{equation}	
and 
\begin{equation}
\sum_{k=0}^{m}x^{k}=\frac{1}{1-x}-\frac{x^{m+1}}{1-x}.\label{40}
\end{equation}
From \eqref{39} and \eqref{40}. we note that 
\begin{align}
\sum_{k=1}^{m}(k)_{n,\lambda}x^{k}&=x^{n\lambda}\bigg(x^{1-\lambda}\frac{d}{dx}\bigg)^{n}\sum_{k=1}^{m}x^{k}\nonumber\\
&=x^{n\lambda}\bigg(x^{1-\lambda}\frac{d}{dx}\bigg)^{n}\frac{1}{1-x}-x^{n\lambda}\bigg(x^{1-\lambda}\frac{d}{dx}\bigg)^{n}\frac{x^{m+1}}{1-x}\label{41}\\
&=\frac{xA_{n,\lambda}(x)}{(1-x)^{n+1}}-x^{m+2}\sum_{l=0}^{n}\binom{n}{l} (m+1)_{n-l,\lambda}\frac{A_{l,\lambda}(x)}{(1-x)^{l+1}}.\nonumber
\end{align}
Therefore, by \eqref{41}, we obtain the following theorem.

 \begin{theorem}
 For $m,n\ge 1$, we have 
 \begin{displaymath}
 	\sum_{k=1}^{m}(k)_{n,\lambda}x^{k}=\frac{xA_{n,\lambda}(x)}{(1-x)^{n+1}}-x^{m+2}\sum_{l=0}^{n}\binom{n}{l}(m+1)_{n-l,\lambda}\frac{A_{l,\lambda}(x)}{(1-x)^{l+1}}.
 \end{displaymath}	
 \end{theorem}
From \eqref{36}, we note that 
\begin{align}
\frac{t-1}{t-e_{-\lambda}\big((t-1)x\big)}&=\frac{t-1}{e_{-\lambda}((t-1)x)\big(te_{\lambda}((1-t)x)-1)}\label{42} \\
&=(1-t)e_{\lambda}((1-t)x)\sum_{j=0}^{\infty}t^{j}e_{\lambda}^{j}((1-t)x)\nonumber \\
&=(1-t)\sum_{j=0}^{\infty}t^{j}e_{\lambda}^{j+1}((1-t)x)\nonumber \\
&=\sum_{n=0}^{\infty}(1-t)^{n+1}\bigg(\sum_{j=0}^{\infty}t^{j}(j+1)_{n,\lambda}\bigg)\frac{x^{n}}{n!}\nonumber \\
&=\sum_{n=0}^{\infty}(1-t)^{n+1}\frac{A_{n,\lambda}(t)}{(1-t)^{n+1}}\frac{x^{n}}{n!}=\sum_{n=0}^{\infty}A_{n,\lambda}(t)\frac{x^{n}}{n!}. \nonumber
\end{align}
Therefore we obtain the following theorem.
\begin{theorem}
The exponential generating function of the degenerate Eulerian polynomials is given by 
\begin{displaymath}
\frac{t-1}{t-e_{-\lambda}((t-1)x)}=\sum_{n=0}^{\infty}A_{n,\lambda}(t)\frac{x^{n}}{n!}.
\end{displaymath}	
\end{theorem}
From Theorem 9, we note that 
\begin{align}
t-1&=\bigg(\sum_{l=0}^{\infty}A_{l,\lambda}(t)\frac{x^{l}}{l!}\bigg)\big(t-e_{-\lambda}((t-1)x)\big)\label{43}\\
&=\sum_{n=0}^{\infty}tA_{n,\lambda}(t)\frac{x^{n}}{n!}-\sum_{n=0}^{\infty}\bigg(\sum_{k=0}^{n}\binom{n}{k}A_{k,\lambda}(t)(t-1)^{n-k}\langle 1\rangle_{n-k,\lambda}\bigg)\frac{x^{n}}{n!} \nonumber 	\\
&=\sum_{n=0}^{\infty}\bigg\{tA_{n,\lambda}(t)-\sum_{k=0}^{n}\binom{n}{k}A_{k,\lambda}(t)(t-1)^{n-k}\langle 1\rangle_{n-k,\lambda}\bigg\}\frac{x^{n}}{n!}. \nonumber
\end{align}
By comparing the coefficients on both sides of \eqref{43}, we get 
\begin{equation}
tA_{n,\lambda}(t)-\sum_{k=0}^{n}\binom{n}{k}A_{k,\lambda}(t)(t-1)^{n-k}\langle 1\rangle_{n-k,\lambda}=(t-1)\delta_{0,n},\label{44}	
\end{equation}
where $\delta_{n,k}$ is the Kronecker's symbol. \par 
Therefore, by \eqref{44}, we obtain the following theorem. 
\begin{theorem}
For $n\in\mathbb{N}$, we have 
\begin{displaymath}
	A_{n,\lambda}(t)=\sum_{k=0}^{n-1}\binom{n}{k}A_{k,\lambda}(t)(t-1)^{n-k-1}\langle 1\rangle_{n-k,\lambda}.
\end{displaymath}	
\end{theorem}
From \eqref{6}, we note that 
\begin{align}
\langle x\rangle_{m+1}&=\sum_{k=0}^{m+1}{m+1 \brack k}_{\lambda}\langle x\rangle_{k,\lambda}= \sum_{k=1}^{m+1}{m+1 \brack k}_{\lambda}\langle x\rangle_{k,\lambda}\label{45} \\
&= \sum_{k=0}^{m}{m+1 \brack k+1}_{\lambda}\langle x\rangle_{k+1,\lambda}= \sum_{k=0}^{m}{m+1 \brack k+1}_{\lambda}\langle x\rangle_{k,\lambda}(x+k\lambda)\nonumber
\end{align}
Thus, by \eqref{45}, we get 
\begin{align}
&\sum_{k=0}^{m}{m+1 \brack k+1}_{\lambda}\langle x\rangle_{k,\lambda}=\frac{1}{x}\langle x\rangle_{m+1}-\frac{\lambda}{x}\sum_{k=0}^{m}{m+1 \brack k+1}_{\lambda}k\langle x\rangle_{k,\lambda}\label{46}\\
&=\langle x+1\rangle_{m}-\frac{\lambda}{x}\sum_{k=1}^{m}{m+1 \brack k+1}_{\lambda}k\langle x\rangle_{k,\lambda}\nonumber \\
&=\langle x+1\rangle_{m}-\frac{\lambda}{x}\sum_{k=0}^{m-1}{m+1 \brack k+2}_{\lambda}(k+1)\langle x\rangle_{k+1,\lambda}\nonumber \\
&=\sum_{k=0}^{m}{m \brack k}_{\lambda}\langle x+1\rangle_{k,\lambda}-\lambda
\sum_{k=0}^{m-1}{m+1 \brack k+2}_{\lambda}(k+1)\langle x+1
\rangle_{k,\lambda}\nonumber \\
&=\sum_{k=0}^{m}{m \brack k}_{\lambda}\sum_{j=0}^{k}\binom{k}{j}\langle 1\rangle_{k-j,\lambda}\langle x\rangle_{j,\lambda}-\lambda\sum_{k=0}^{m}{m+1 \brack k+2}_{\lambda}(k+1)\sum_{j=0}^{k}\binom{k}{j}\langle 1\rangle_{k-j,\lambda}\langle x\rangle_{j,\lambda}\nonumber\\
&=\sum_{j=0}^{m}\bigg\{\sum_{k=j}^{m}\binom{k}{j}\bigg({m \brack k}_{\lambda}\langle 1\rangle_{k-j,\lambda}-\lambda(k+1){m+1\brack k+2}_{\lambda}\langle 1\rangle_{k-j,\lambda}\bigg)\langle x\rangle_{j,\lambda}.\nonumber
\end{align}
Comparing the coefficients on both sides of \eqref{46}, we obtain the following theorem. 
\begin{theorem}
For $j,m\in\mathbb{Z}$ with $j,m\ge 0$, we have 
\begin{displaymath}
	{m+1 \brack j+1}_{\lambda}=\sum_{k=j}^{m}\binom{k}{j}\bigg\{{m \brack k}_{\lambda}\langle 1\rangle_{k-j,\lambda}-\lambda(k+1){m+1 \brack k+2}_{\lambda}\langle 1\rangle_{k-j,\lambda}\bigg\}.
\end{displaymath}	
\end{theorem}

\section{Conclusion}
Eulerian numbers and polynomials have long been studied due to their number-theoretic and combinatorial importance.
In this paper, we considered a finite sum involving the generalized falling factorial polynomials, in some special cases of which we represented it in terms of the degenerate Stirling numbers of the second kind, the degenerate Bernoulli polynomials and the degenerate Frobenius-Euler polynomials. Also, we considered the degenerate Eulerian polynomials and deduced the generating function and a recurrence relation for them. \par
There are various ways of studying special numbers and polynomials, to mention a few, generating functions, $p$-adic analysis, umbral calculus, combinatorial methods, differential equations, special functions, probability theory and analytic number theory. It is one of our future projects to continue to explore various degenerate versions of many special polynomials and numbers.

\end{document}